\documentclass[12pt]{article}
\usepackage{amsfonts}
\usepackage{bbm}
\usepackage{mathrsfs}
\usepackage{amssymb}

\def\hang{\hangindent\parindent}
\def\textindent#1{\indent\llap{#1\enspace}\ignorespaces}

\title{Valuation Extensions of Algebras \\ Defined by Monic  Gr\"obner Bases
\thanks{Project supported by the National Natural Science Foundation
of China (10971044).}}

\author{Huishi Li\\
{\small Department of Applied Mathematics}\\
{\small College of Information Science and Technology}\\
{\small Hainan University}\\
{\small  Haikou 570228, China}}

\date{}

\begin{document}
\maketitle
\begin{center}
\begin{minipage}{120mm}
{\small {\bf Abstract.} Let $K$ be a field, $\mathcal {O}_v$ a
valuation ring of $K$ associated to a valuation $v$:
$K\rightarrow\Gamma\cup\{\infty\}$, and ${\bf m}_v$ the unique
maximal ideal of $\mathcal {O}_v$. Consider an ideal $\mathcal {I}$
of the free $K$-algebra $K\langle X\rangle =K\langle
X_1,...,X_n\rangle$ on $X_1,...,X_n$. If  ${\cal I}$ is generated by
a subset $\mathcal {G}\subset{\cal O}_v\langle X\rangle$ which is a
monic Gr\"obner basis of ${\cal I}$ in $K\langle X\rangle$, where
$\mathcal {O}_v\langle X\rangle =\mathcal{O}_v\langle
X_1,...,X_n\rangle$ is the free $\mathcal{O}_v$-algebra on
$X_1,...,X_n$, then the valuation $v$ induces naturally an
exhaustive and separated $\Gamma$-filtration $F^vA$ for the
$K$-algebra $A=K\langle X\rangle /\mathcal {I}$, and moreover
$\mathcal{I}\cap\mathcal{O}_v\langle X\rangle
=\langle\mathcal{G}\rangle$ holds in $\mathcal{O}_v\langle
X\rangle$; it follows that, if furthermore $\mathcal{G}\not\subset
{\bf m}_v{O}_v\langle X\rangle$ and $k\langle X\rangle
/\langle\overline{\mathcal G}\rangle$ is a domain, where
$k=\mathcal{O}_v/{\bf m}_v$ is the residue field of $\mathcal{O}_v$,
$k\langle X\rangle =k\langle X_1,...,X_n\rangle$ is the free
$k$-algebra on $X_1,...,X_n$, and $\overline{\mathcal G}$ is the
image of $\mathcal{G}$ under the canonical epimorphism
$\mathcal{O}_v\langle X\rangle\rightarrow k\langle X\rangle$, then
$F^vA$ determines a valuation function $A\rightarrow
\Gamma\cup\{\infty\}$, and thereby  $v$ extends naturally to a
valuation function on the (skew-)field $\Delta$ of fractions of $A$
provided $\Delta$ exists. }
\end{minipage}\end{center} {\parindent=0pt\vskip 6pt

{\bf 2000 Mathematics Classification} Primary 16W60; Secondary 16Z05
(68W30).\vskip 6pt

{\bf Key words} Filtered algebra, graded algebra, Gr\"obner basis,
valuation.}

\def\QED{\hfill{$\Box$}}
\def \r{\rightarrow}\def\NZ{\mathbb{N}}

\def\mapright#1#2{\smash{\mathop{\longrightarrow}\limits^{#1}_{#2}}}

\def\v5{\vskip .5truecm}

\def\OV#1{\overline {#1}}
\def\item{\par\hang\textindent}

\def\LH{{\bf LH}}\def\LM{{\bf LM}}\def\LT{{\bf
LT}}\def\KX{K\langle X\rangle} \def\KS{K\langle X\rangle}
\def\B{{\cal B}} \def\LC{{\bf LC}} \def\G{{\cal G}} \def\FRAC#1#2{\displaystyle{\frac{#1}{#2}}}
\def\SUM^#1_#2{\displaystyle{\sum^{#1}_{#2}}} \def\O{{\cal O}}  \def\I{{\cal I}} \def\DX{D\langle X\rangle}
\def\M{{\bf m}} \def\RS{R\langle X\rangle}

\def\KS{K\langle X\rangle}
\def\LR{\langle X\rangle}
\def\HL{{\rm LH}}
\vskip 1truecm

\section*{1. Introduction}
In the so-called noncommutative algebraic geometry, the class of
schematic algebras in the sense of ([VOW1, [VOW3]) has provided an
ample stage to play on. Among others of the topics concerning
noncommutative geometric objects associated to schematic algebras,
noncommutative valuations are applied to obtain tools for an
equivalent of divisor theory in noncommutative geometry (the reader
is referred to [VO] for details on this aspect). In the study of
extending commutative valuations to noncommutative valuations,
filtered-graded structural methods have been used successfully to
obtain sufficient conditions assuring the existence of an extension
([LVO2], [MVO], [Li1], [VOW2], [VO], [BVO]). More precisely, let $K$
be a field and $\O_v$ a  valuation ring of $K$ associated to a
valuation $v$: $K\r\Gamma\cup\{\infty\}$, where $\Gamma$ is a
totally ordered abelian additive group. Then $v$ determines an
exhaustive and separated $\Gamma$-filtration $F^vK=\{
F^v_{\gamma}K\}_{\gamma\in\Gamma}$ for $K$, where $F^v_{\gamma}K=\{
\lambda~|~\lambda\in K,~v(\lambda)\ge -\gamma )\}$, such that
$F^v_{\gamma_1}K\cdot F^v_{\gamma_2}K=F^v_{\gamma_1+\gamma_2}K$ for
all $\gamma_1,\gamma_2\in\Gamma$, i.e., $F^vK$ is a strong
$\Gamma$-filtration in the sense of [LVO1]. Consider an affine
$K$-algebra $A=K[a_1,...,a_n]$ with the (finite or infinite) set of
defining relations $\G$, that is, $A\cong \KS /\I$ with $\I =\langle
\G\rangle$, where $\KS =K\langle X_1,...,X_n\rangle$ is the free
$K$-algebra on $X_1,...,X_n$. From loc. cit. we have learnt that the
key points  of naturally extending the given valuation $v$ of $K$ to
$A$ and further to the (skew-)field $\Delta$ of fractions of $A$
(provided $\Delta$ exists) are to assure
that{\parindent=.8truecm\par

\item{(1)} the valuation $\Gamma$-filtration $F^vK$ of $K$ induces an exhaustive $\Gamma$-filtration
$F^vA=\{ F^v_{\gamma}A\}_{\gamma\in\Gamma}$ for $A$ in a natural
way, i.e., $F^v_{\gamma}K=K\cap F^v_{\gamma}A$ for every
$\gamma\in\Gamma$, such that $ F_0^vA=\O_v\langle X\rangle +\I /\I$
and $F^v_{<0}A=\M_vF^v_0A$, where $\O_v\langle X\rangle =\O_v\langle
X_1,...,X_n\rangle$ is the free $\O_v$-algebra on $X_1,...,X_n$, and
$\M_v$ is the unique maximal ideal of $\O_v$;\par

\item{(2)} the $\Gamma$-filtration $F^vA$ obtained in (1) above is separated, i.e.,
$0\ne a\in A$ implies that there is some $\Gamma\in\Gamma$ such that
$a\in F^v_{\gamma}A-F^v_{<\gamma}A$, in particular $1\in
F^v_0A-F^v_{<0}A$,  where $F^v_{<\gamma}A=\cup_{\gamma
'<\gamma}F^v_{\gamma '}A$; and
\par

\item{(3)} if $\G\subset \O_v\langle X\rangle$ then $\O_v\langle X\rangle\cap
\I =\langle \G\rangle$ holds in $\O_v\langle X\rangle$. In loc. cit.
this property is referred to as saying that the $\O_v$-algebra
$\O_v\langle X\rangle +\I /\I$ defines a good reduction for the
$K$-algebra $A$. {\parindent=0pt\par

For a connected positively $\NZ$-graded $K$-algebra $A$, it was
shown in ([VO], Theorem 4.3.7; [BVO], Theorem 2.2) that the
$\Gamma$-filtration $F^vA$ constructed in loc. cit. may have the
properties (1) -- (2) provided $A$ has a PBW $K$-bsis in the
classical sense; while the property (3) may be derived under the
so-called $v$-comaximal condition assumed on the ideal
$I=\langle\G\rangle$ of $\O_v\langle X\rangle$, i.e., $I\cap
(F^v_{\gamma}K)\langle X\rangle =(F^v_{\gamma}K)I$ for every
$\gamma\in\Gamma$ ([MVO], Lemma 2.1;  [VO], Lemmaa 4.3.2). }}\par

From ([Li2], CH.III Theorem 1.5; [Li3], Theorem 3.1) we know that,
for algebras of the type $A=\KS /\I$ as considered above, the
property that $A$ has a classical PBW $K$-basis may be equivalent to
the property that $\I$ is generated by a (finite or infinite)
Gr\"obner basis of special type. For instance, all the concrete
algebras quoted in ([MVO], [Li1], [VO], [BVO]) are indeed defined by
Gr\"obner bases that give rise to PBW $K$-bases (cf. [Li2], [Li3]).
Inspired by such a fact, we aim to demonstrate the following main
result in this paper:{\parindent=.5truecm\vskip 6pt

\item{$\bullet$} If $\G\subset \O_v\langle X\rangle$ forms a monic Gr\"obner basis for the ideal $\I$ in
$\KS$, where ``monic" means that the leading coefficient of every
element in $\mathcal{G}$ is 1 (see Section 2 for details), then $A$
has the three properties (1) -- (3) described above. It follows
that, if furthermore $\G\not\subset\M_v\O_v\langle X\rangle$ and the
$k$-algebra $k\langle X\rangle /\langle\OV{\G}\rangle$ is a domain,
where $k=\O_v/\M_v$ is the residue field of $\O_v$, $k\langle
X\rangle =k\langle X_1,...,X_n\rangle$ is the free $k$-algebra on
$X_1,...,X_n$, and $\OV{\G}$ is the canonical image of $\G$ in
$\O_v\langle X\rangle /\M_v\O_v\langle X\rangle$, then $F^vA$
determines a valuation function $A\rightarrow \Gamma\cup\{\infty\}$,
and thereby  $v$ extends naturally to a valuation function on the
(skew-)field $\Delta$ of fractions of $A$ provided $\Delta$ exists.
{\parindent=0pt\vskip 6pt

The result mentioned above will be reached by deriving several
results for $R$-algebras over an arbitrary commutative ring $R$,
where the filtration considered will be $\Gamma$-filtration with
$\Gamma$ a totally ordered (commutative or noncommutative) monoid.
That is, the results obtained in Sections 3 -- 5 may be of
independent interest, for instance, they may be used to study
valuation extensions of commutative algebras defined by monic
Gr\"obner bases (see the remark given at the end of this paper), and
they may also be used to study more general reductions of algebras
over a field $K$ as specified in [LVO3].
\parindent=.5truecm\par

By the algorithmic Gr\"obner basis theory for free $K$-algebras over
a field $K$ ([Mor], [Gr]), in principle every finitely presented
algebra $A=\KS /\I$ has the defining ideal $\I$ generated by a
(finite or infinite) Gr\"obner basis $\G$ which can always be
assumed to be monic. Furthermore, by [Li3] (or see Proposition 2.7
in Section 2 below), if $D$ is a subring of $K$ with the  same
multiplicative identity 1, then $\G\subset D\langle X\rangle
=D\langle X_1,...,X_n\rangle$ is a monic Gr\"obner basis for the
ideal $I=\langle\G\rangle$ in $D\langle X\rangle$ if and only if
$\G$ is a monic Gr\"obner basis for the ideal $\I =\langle
\G\rangle$ in $\KS$ with respect to the same monomial ordering on
both $\KS$ and $D\langle X\rangle$. In this sense, the work of this
paper may be viewed as a computational approach to solving the
valuation extension problem. So, the contents of this paper are
organized as follows.\par

1. Introduction\par

2. Monic Gr\"obner bases over rings\par

3.  Extending $FR$ naturally to $FA$  by Gr\"obner bases over
$F_0R$\par

4. Realizing the separability of $FA$ by Gr\"obner bases over
$F_0R$\par

5. Realizing good reductions for $A$ by Gr\"obner bases over
$D\subset R$\par

6. Realizing valuation extensions of $A$ by Gr\"obner bases over
$\O_v${\parindent=0pt\par

Unless otherwise stated, rings considered in this paper are
associative rings with multiplicative identity 1,  ideals are meant
two-sided ideals, and modules are unitary left modules. For a subset
$U$ of a ring $S$, we write $\langle U\rangle$  for the ideal
generated by $U$. Moreover, we use $\NZ$, respectively $\mathbb{Z}$,
to denote the set of nonnegative integers, respectively the set of
integers. Moreover, valuations of a (skew-)field $\Delta$ are in the
sense of O. Schilling [Sc]}.

\section*{2. Monic Gr\"obner Bases over Rings}
For the reader's convenience, in this section we briefly recall from
[Li3] some basics on monic Gr\"obner bases in free algebras over
rings. Classical Gr\"obner basis theory for free algebras over a
field $K$ is referred to [Mor] and [Gr]. \v5

Let $R$ be an arbitrary commutative ring, $\RS =R\langle
X_1,...,X_n\rangle$ the free $R$-algebra of $n$ generators, and
$\B_R$ the standard $R$-basis of $\RS$ consisting of monomials
(words in alphabet $X=\{ X_1,...,X_n\}$, including empty word  which
is identified with the multiplicative identity element 1 of $\RS$).
Unless otherwise stated, monomials in $\B_R$ are denoted by lower
case letters $u,v,w,s,t,\cdots$. By a {\it monomial ordering} on
$\B_R$ (or on $\RS$) we mean a well-ordering $\prec$ on $\B_R$ which
satisfies: \par

(M1) For $w,u,v,s\in\B_R$, $u\prec v$ implies $wus\prec wvs$;\par

(M2) For $w,u,v\in\B_R$, $w=uv$ implies $u\preceq w$ and $v\preceq
w$. {\parindent=0pt\par

In particular, by an $\NZ$-{\it graded monomial ordering} on $\B_R$,
denoted $\prec_{gr}$, we mean a monomial ordering on $\B_R$ which is
defined subject to a well-ordering $\prec$ on $\B_R$, that is,  for
$u,v\in\B_R$, $u\prec_{gr} v$ if either deg$(u)<$ deg$(v)$ or
deg$(u)=$ deg$(v)$ but $u\prec v$, where deg$(~)$ denotes the degree
function on elements of $\RS$ with respect to a fixed {\it weight}
$\NZ$-{\it gradation} of $\RS$ (i.e. each $X_i$ is assigned  a
positive degree $n_i$, $1\le i\le n$). For instance, the usual
$\NZ$-graded (reverse) lexicographic ordering is a popularly used
$\NZ$-graded monomial ordering.}\par

If $\prec$ is a monomial ordering on $\B_R$ and
$f=\sum_{i=1}^s\lambda_iw_i\in \RS$, where $\lambda_i\in R-\{ 0\}$
and $w_i\in\B_R$, such that $w_1\prec w_2\prec\cdots\prec w_s$, then
the {\it leading monomial} of $f$ is defined as $\LM (f)=w_s$ and
the {\it leading coefficient} of $f$ is defined as $\LC
(f)=\lambda_s$. For a subset $H\subset\RS$, we write $\LM (H)=\{\LM
(f)~|~f\in H\}$ for the set of leading monomials of $S$. We say that
a subset $G\subset\RS$ is {\it monic} if $\LC (g)=1$ for every $g\in
G$. Moreover, for $u,v\in\B_R$, we say that $v$ divides $u$, denoted
$v|u$, if $u=wvs$ for some $w,~s\in\B_R$. \par

With notation and all definitions as above, it is easy to see that a
division algorithm by a monic subset $G$ is valid in $\RS$ with
respect to any fixed monomial ordering $\prec$ on $\B_R$. More
precisely, let $f\in\RS$. Noticing $\LC (g)=1$ for every $g\in G$,
if $\LM (g)|\LM (f)$ for some $g\in G$, then $f$ can be written as
$f=\LC (f)ugv+f_1$ with $u$, $v\in\B_R$, $f_1\in\RS$ satisfying $\LM
(f_1)\prec\LM (f)$; if $\LM (g){\not |}~\LM (f)$ for all $g\in G$,
then $f=f_1+\LC (f)\LM (f)$ with $f_1=f-\LC (f)\LM (f)$ satisfying
$\LM (f_1)\prec\LM (f)$. Next, consider the divisibility of $\LM
(f_1)$ by $\LM (g)$ with $g\in G$, and so forth. Since $\prec$ is a
well-ordering, after a finite number of successive division by
elements in $G$ in this way, we see that $f$ can be written as
$$\begin{array}{rcl} f&=&\sum_{i,j}\lambda_{ij}u_{ij}g_jv_{ij}+r_f,~\hbox{where}~\lambda_{ij}\in R,~u_{ij},v_{ij}\in\B_R,~g_j\in G,\\
&{~}&\hbox{and}~r_f=\sum_p\lambda_pw_p~\hbox{with}~\lambda_p\in
R,~w_p\in\B_R,\\
&{~}&\hbox{satisfying}~\LM (u_{ij}g_jv_{ij})\preceq\LM
(f)~\hbox{whenever}~\lambda_{ij}\ne 0,\\
&{~}&\LM (r_f)\preceq\LM (f)~\hbox{and}~\LM (g)\not |~w_p~\hbox{for
every}~g\in G~\hbox{whenever}~\lambda_p\ne 0.\end{array}$$
If, $r_f=0$ in the representation of $f$ obtained above, then we say
that $f$ {\it is reduced to} 0 {\it by division by} $G$, and we
write $\OV f^G=0$ for this property. The validity of such a division
algorithm by $G$ leads to the following definition.
{\parindent=0pt\v5

{\bf 2.1. Definition} Let $\prec$ be a fixed monomial ordering on
$\B_R$, and $I$ an ideal of $\RS$. A {\it monic Gr\"obner basis} of
$I$ is a subset $\G\subset I$ satisfying:\par

(1) $\G$ is monic; and
\par

(2) $f\in I$ and $f\ne 0$ implies $\LM (g)|\LM (f)$ for some $g\in
\G$. }\v5

By the division algorithm presented above, it is clear that a monic
Gr\"obner basis of $I$ is first of all a generating set of the ideal
$I$, i.e., $I=\langle\G\rangle$, and moreover, a monic Gr\"obner
basis of $I$ can be characterized as follows. {\parindent=0pt\v5

{\bf 2.2. Proposition}  Let $\prec$ be a fixed monomial ordering on
$\B_R$, and $I$ an ideal of $\RS$. For a monic subset $\G\subset I$,
the following statements are equivalent:\par (i) $\G$ is a monic
Gr\"obner basis of $I$;\par (ii) Each nonzero $f\in I$ has a
Gr\"obner representation:
$$\begin{array}{rcl} f&=&\sum_{i,j}\lambda_{ij}u_{ij}g_jv_{ij},
~\hbox{where}~\lambda_{ij}\in R,~u_{ij},v_{ij}\in\B_R,~g_j\in G,\\
&{~}&\hbox{satisfying}~\LM (u_{ij}g_jv_{ij})\preceq\LM
(f)~\hbox{whenever}~\lambda_{ij}\ne 0,\end{array}$$ or equivalently,
$\OV f^{\G}=0$; \par (iii) $\langle\LM (\G )\rangle =\langle\LM
(I)\rangle$. \QED}\v5

Let $\prec$ be a monomial ordering on the standard $R$-basis $\B_R$
of $\RS$, and let $G$ be a monic subset of $\RS$. We call an element
$f\in \RS$ a {\it normal element} (mod $G$) if $f=\sum_j\mu_jv_j$
with $\mu_j\in R$, $v_j\in\B_R$, and $f$ has the property that $\LM
(g){\not |}~v_j$ for every $g\in G$ and every $\mu_j\ne 0$. The set
of normal monomials in $\B_R$ (mod $G$) is denoted by $N(G)$, i.e.,
$$N(G)=\{ u\in\B_R~|~\LM (g){\not |}~u, ~g\in G\} .$$
Thus, an element $f\in\RS$ is normal (mod $G$) if and only if $f\in
\sum_{u\in N(G)}Ru$. {\parindent=0pt\v5

{\bf 2.3. Proposition} Let $\G$ be a monic Gr\"obner basis of the
ideal $I=\langle\G\rangle$ in $\RS$ with respect to some monomial
ordering $\prec$ on $\B_R$. Then each nonzero $f\in \RS$ has a
finite presentation
$$f=\sum_{i,j}\lambda_{ij}s_{ij}g_iw_{ij}+r_f,\quad \lambda_{ij}\in R,~s_{ij},
w_{ij}\in\B_R,~g_i\in \G,$$ where $\LM (s_{ij}g_iw_{ij})\preceq\LM
(f)$ whenever $\lambda_{ij}\ne 0$, and either $r_f=0$ or $r_f$ is a
unique normal element (mod $\G$). Hence, $f\in I$ if and only if
$r_f=0$, solving the ``membership problem" for $I$.\par \QED} \v5

The foregoing results enable us to obtain further characterization
of a monic Gr\"obner basis $\G$, which, in turn, gives rise to the
fundamental decomposition theorem of the $R$-module $\RS$ by the
ideal $I=\langle\G\rangle$, and thereby yields a free $R$-basis for
the $R$-algebra $\RS /I$. {\parindent=0pt\v5

{\bf 2.4. Theorem}  Let  $I=\langle\G\rangle$ be an ideal of $\RS$
generated by a monic subset $\G$. With notation as above, the
following statements are equivalent.\par (i) $\G$ is a monic
Gr\"obner basis of $I$. \par (ii) The $R$-module $\RS$ has the
decomposition
$$\RS=I\oplus \sum_{u\in N(\G )}Ru=\langle\LM (I)\rangle\oplus\sum_{u\in N(\G )}Ru.$$
(iii) The canonical image $\OV{N(\G )}$ of $N(\G )$ in $\RS
/\langle\LM (I)\rangle$ and $\RS/I$ forms a free $R$-basis for $\RS
/\langle\LM (I)\rangle$ and $\RS/I$ respectively. \par\QED}\v5

Before mentioning a version of the termination theorem in the sense
of ([Mor], [Gr]) for verifying an LM-reduced monic Gr\"obner basis
in $\RS$ (see the definition below), we need a little more
preparation. \v5

Given a monomial ordering $\prec$ on $\B_R$, we say that a subset
$G\subset\RS$ is {\it LM-reduced} if
$$\LM (g_i)\not |~\LM (g_j)~\hbox{for all}~g_i,g_j\in G~\hbox{with}~g_i\ne g_j.$$
If a subset $G\subset\RS$ is both LM-reduced and monic, then we call
$G$ an {\it LM-reduced monic subset}. Thus we have the notion of an
{\it LM}-{\it reduced monic Gr\"obner basis}.\par Let $I$ be an
ideal of $\RS$. If $\G$ is a monic Gr\"obner basis of $I$ and
$g_1,g_2\in\G$ such that $g_1\ne g_2$ but $\LM (g_1)|\LM (g_2)$,
then clearly $g_2$ can be removed from $\G$ and the remained subset
$\G-\{ g_2\}$ is again a monic Gr\"obner basis for $I$. Hence, in
order to have a better criterion for monic Gr\"obner basis we need
only to consider the subset which is both LM-reduced and monic. \v5

Let $\prec$ be a monomial ordering on $\B_R$. For two monic elements
$f,~g\in\RS-\{ 0\}$, including $f=g$, if there are monomials
$u,v\in\B_R$ such that \par

(1) $\LM (f)u=v\LM (g)$, and\par

(2) $\LM (f)\not |~v$ and $\LM (g)\not |~u$,{\parindent=0pt\par then
the element
$$o(f,u;~v,g)=f\cdot u-
v\cdot g$$ is called an {\it overlap element} of $f$ and $g$. From
the definition it is clear that
$$\LM ((o(f,u;~v,g))\prec \LM (fu)=\LM (vg),$$
and moreover, there are only finitely many overlap elements for each
pair $(f,g)$ of monic elements in $\RS$. So, for a finite subset of
monic elements $\G\subset\RS$, actually as in the classical case
([Mor], [Gr]), the termination theorem below enables us to check
effectively whether $\G$ is a Gr\"obner basis of $I$ or not.\v5

{\bf 2.5. Theorem} (Termination theorem)  Let $\prec$ be a fixed
monomial ordering on $\B_R$. If $\G$ is an LM-reduced monic subset
of $\RS$, then $\G$ is an LM-reduced monic Gr\"obner basis for the
ideal $I=\langle\G\rangle$ if and only if for each pair
$g_i,g_j\in\G$, including $g_i=g_j$, every overlap element
$o(g_i,u;~v,g_j)$ of $g_i$, $g_j$ has the property
$\overline{o(g_i,u;~v,g_j)}^{\G}=0,$ that is, by division by $\G$,
every $o(g_i,u;~v,g_j)$ is reduced to zero.\par \QED\v5

{\bf Remark} (i) Obviously, if $\G\subset\RS$ is an LM-reduced
subset with the property that each $g\in \G$ has the leading
coefficient $\LC (g)$ which is invertible in $R$, then Theorem 2.5
is also valid for $\G$.\par

(ii) It is obvious as well that Theorem 2.5 does not necessarily
induce an analogue of the Buchberger algorithm as in the classical
case. \par

(iii) It is not difficult to see that all results we presented so
far are valid for getting monic Gr\"obner bases in a commutative
polynomial ring $R[x_1,...,x_n]$ over an arbitrary commutative ring
$R$ where overlap elements are replaced by $S$-polynomials.}\v5

By virtue of Theorem 2.5 (or more precisely, its proof given in
[Li3]), the following two propositions are
obtained.{\parindent=0pt\v5

{\bf 2.6. Proposition} Let $K\langle X\rangle =K\langle
X_1,...,X_n\rangle$ be the free algebra of $n$ generators over a
field $K$, and let $R\langle X\rangle =R\langle X_1,...,X_n\rangle$
be the free algebra of  $n$ generators over an arbitrary commutative
ring $R$. With notation as before, fixing the same monomial ordering
$\prec$ on both $\KS$ and $\RS$, the following statements hold.\par
(i) If a monic subset $\G\subset\KS$ is a Gr\"obner basis for the
ideal $\langle\G\rangle$ in $\KS$, then, taking a counterpart of
$\G$ in $\RS$ (if it exists), again denoted by $\G$, $\G$ is a monic
Gr\"obner basis for the ideal $\langle\G\rangle$ in $\RS$. \par

(ii) If a monic subset $\G\subset\RS$ is a Gr\"obner basis for the
ideal $\langle\G\rangle$ in $\RS$, then, taking a counterpart of
$\G$ in $\KS$ (if it exists), again denoted by $\G$, $\G$ is a
Gr\"obner basis for the ideal $\langle\G\rangle$ in $\KS$.\par
\QED\v5

{\bf 2.7. Proposition} Let $R$ be a commutative ring and $R'$ a
subring of $R$ with the same identity element 1. If we consider the
free $R$-algebra $\RS =R\langle X_1,...,X_n\rangle$ and the free
$R'$-algebra $R'\langle X\rangle =R'\langle X_1,...,X_n\rangle$,
then the following two statements are equivalent for a subset
$\G\subset R'\langle X\rangle$:\par

(i) $\G$ is an LM-reduced monic Gr\"obner basis for the ideal
$I=\langle\G\rangle$ in $R'\langle X\rangle$ with respect to some
monomial ordering $\prec$ on the standard $R'$-basis $\B_{R'}$ of
$R'\langle X\rangle$;\par

(ii) $\G$ is an LM-reduced monic Gr\"obner basis for the ideal
$J=\langle\G\rangle$ in $\RS$ with respect to the monomial ordering
$\prec$ on the standard $R$-basis $\B_R$ of $\RS$, where $\prec$ is
the same monomial ordering used in (i).\par\QED }\v5

Let $K$ be a field. From the literature we know that numerous
well-known $K$-algebras, such as  Weyl algebras over $K$, enveloping
algebras of $K$-Lie algebras, exterior $K$-algebras, Clifford
$K$-algebras, down-up $K$-algebras, quantum binomial $K$-algebras,
most popularly studied quantum groups over $K$, etc., all have
defining relations that form an LM-reduced monic  Gr\"obner basis in
free $K$-algebras (cf. [Li2], [Laf], [G-I]). Hence, by Proposition
2.6, if the field $K$ is replaced by a commutative ring $R$, then
all of these $R$-algebras (if they exist) have defining relations
that form an LM-reduced monic Gr\"obner basis in a free $R$-algebra.
The reader is referred to [Li3] for more details on this topic and
for more concrete examples. \v5

\section*{3. Extending $FR$ Naturally to $FA$ by Gr\"obner Bases over $F_0R$}
Let $R$ be an arbitrary {\it commutative ring}, and $\Gamma$ a
totally ordered (commutative or noncommutative) {\it monoid} with
the total ordering $<$. To make the notation uniform in our context,
we first fix the convention:{\parindent=.5truecm\par

\item{$\bullet$} From now on in this paper we use $+$ to denote the binary
operation of $\Gamma$, and we use $0$ to denote the neutral element
of $\Gamma$ (though $\Gamma$ is not necessarily commutative).

\item{$\bullet$} The definition of an exhaustive $\Gamma$-filtration defined for $R$ below
applies to every $R$-algebra (ring) considered in this paper.
\parindent=0pt\par

We say that $R$ is equipped with an exhaustive  $\Gamma$-filtration
$FR=\{ F_{\gamma}R\}_{\gamma\in\Gamma}$, if each $F_{\gamma}R$ is an
additive subgroup of $R$ and $FR$ satisfies}\par

(F1) $R=\cup_{\gamma\in\Gamma}F_{\gamma}R$;\par

(F2) $F_{\gamma_1}R\subseteq F_{\gamma_2}R$ whenever
$\gamma_1<\gamma_2$;\par

(F3) $F_{\gamma}R\cdot F_{\tau}R\subseteq F_{\gamma+\tau}R$ for all
$\gamma ,\tau\in\Gamma$;\par

(F4) $1\in F_0R$. {\parindent=0pt\par

Note that $F_0R$ is a subring of $R$ with the same identity element
1. To simplify notation, we write $R_0$ for $F_0R$. Let $\RS
=R\langle X_1,...,X_n\rangle$ be the free $R$-algebra on
$X_1,...,X_n$, and $\I$ an ideal of $\RS$. Considering the
$R$-algebra $A=\RS /\I$, if $\G$ is a monic Gr\"obner basis of $\I$
in $\RS$ with respect to a monomial ordering on the standard
$R$-basis $B_R$ of $\RS$, and if $N(\G )$ denotes the set of normal
monomials in $\B_R$ (mod $\G$) (see Section 2), then, by Theorem
2.4, the canonical image $\OV{N(\G )}$ of $N(\G )$ in $A=\RS
/\langle\I \rangle =\RS /\langle\G\rangle$ forms a free $R$-basis
for $A$. Bearing this preliminary in mind, we are able to establish
the following result.\v5

{\bf 3.1. Theorem} Let the commutative ring $R$ be equipped with an
exhaustive  $\Gamma$-filtration $FR=\{
F_{\gamma}R\}_{\gamma\in\Gamma}$. With notation as fixed above,
suppose that the ideal $\I$ is generated by a subset
$\mathcal{G}\subset R_0\langle X\rangle$ which is a monic Gr\"obner
basis of $\I$ in $\RS$ with respect to a monomial ordering $\prec$
on the standard $R$-basis $\B_R$ of $\RS$, where $R_0\langle
X\rangle =R_0\langle X_1,...,X_n\rangle$ is the free $R_0$-algebra
on $X_1,...,X_n$, and without loss of generality we assume that $\LM
(g)\ne 1$ for every $g\in\G$. Then the $R$-algebra $A=\RS /\I=\RS
/\langle\G\rangle$ can be endowed with the exhaustive
$\Gamma$-filtration  $FA=\{ F_{\gamma}A\}_{\gamma\in\Gamma}$ by
putting
$$F_{\gamma}A=\left\{\left. a=\sum_i\lambda_i\bar w_i~\right |~\lambda_i\in F_{\gamma}R,~
\bar w_i\in \OV{N(\G)}\right\} , ~\gamma\in\Gamma ,$$ such that
$F_{\gamma}R=R\cap F_{\gamma}A,~\gamma\in\Gamma,$ that is, $FR$
extends naturally to $FA$. \vskip 6pt

{\bf Proof} We show that $FA$ satisfies the conditions (F1) -- (F4)
required by an exhaustive $\Gamma$-filtration. Obviously each
$F_{\gamma}A$ is an additive subgroup of $A$. Also it is clear that
$F_{\gamma_1}A\subseteq F_{\gamma_2}A$ whenever $\gamma_1<\gamma_2$.
If $0\ne a\in A$, then $a$ can be written uniquely as
$a=\sum_{j=1}^m\mu_j\bar w_j$ with $\mu_j\in R$, $\bar
w_j\in\OV{N(\G )}$. Since $R=\cup_{\gamma\in\Gamma}F_{\gamma}R$, we
may assume that  $\mu_j\in F_{\gamma_j}R$, $1\le j\le m$, and that
$\gamma_1\ge\gamma_2\ge\cdots\ge\gamma_m$ in the totally ordered
monoid $\Gamma$. It follows that $\mu_j\in F_{\gamma_1}R$, $1\le
j\le m$, and thereby $a\in F_{\gamma_1}A$. This shows that
$A=\cup_{\gamma\in\Gamma}F^v_{\gamma}A$. If $a=\sum_i\lambda_i\bar
w_i\in F_{\gamma}A$, $b=\sum_j\mu_j\bar w_j\in F_{\rho}A$, then
$ab=\sum_{i,j}\lambda_i\mu_j\bar w_i\bar w_j$ with
$\lambda_i\mu_j\in F^v_{\gamma +\rho}R$. Since the monic Gr\"obner
basis $\G$ of $\I$ is contained in $R_0\langle X\rangle
=(F_0R)\langle X\rangle$, if we run the the division algorithm of
each $w_iw_j$ by $\G$ in $\RS$, it is indeed implemented in
$R_0\langle X\rangle$. It turns out that
$$w_iw_j=\sum_{p,q}\xi_{pq}u_{pq}g_qv_{pq}+\sum_m\eta_ms_m,~\hbox{where}~\xi_{pq},
\eta_m\in R_0=F_0R,~u_{pq},v_{pq}\in\B_R,~ s_m\in N(\G ).$$
Considering the residue classes in $\RS /\I$, we have $\bar w_i\bar
w_j=\sum_m\eta_m\bar s_m$ with $\eta_m\in F_0R$ and $\bar
s_m\in\OV{N(\G )}$, which implies $ab\in F^v_{\gamma +\rho}A$.
Thereby $F_{\gamma}A\cdot F_{\rho}A\subseteq F_{\gamma+\rho}A$ for
all $\gamma,\rho\in\Gamma$. Moreover, since $1\in N(\G )$ by our
assumption on $\G$, it is easy to see that $1\in F_0A$. This shows
that $FA$ defines an exhaustive $\Gamma$-filtration for
$A$.{\parindent=.5truecm \par

Finally, noticing the fact that $1\in N(\G )$, it is straightforward
that $F_{\gamma}R\subseteq R\cap F_{\gamma}A\subseteq F_{\gamma}R$,
i.e., $F_{\gamma}R=R\cap F_{\gamma}A$ for every $\gamma\in\Gamma$,
as desired. \QED}}\v5

\section*{4. Realizing the Separability of $FA$ by Gr\"obner Bases over $F_0R$}
Let $R$ be an arbitrary commutative ring, and $\Gamma$ a totally
ordered (commutative or noncommutative) monoid with the total
ordering $<$. Suppose that $R$ is equipped with an exhaustive
$\Gamma$-filtration $FR=\{ F_{\gamma}R\}_{\gamma\in\Gamma}$ in the
sense of Section 3. We say that the $\Gamma$-filtration $FR$ is {\it
separated}, if $0\ne \lambda\in R$ implies that there is a
$\gamma\in\Gamma$ such that
$$\lambda\in F_{\gamma}R-F_{<\gamma}R,~\hbox{where}~F_{<\gamma}R=\cup_{\gamma '<\gamma}F_{\gamma'}R.$$
Let $\RS =R\langle X_1,...,X_n\rangle$ be the free $R$-algebra on
$X_1,...,X_n$, $\I$ an ideal of $\RS$, and $A=\RS /\I$. With every
definition and all notations as in Section 3, especially with
$R_0=F_0R$, in this section we aim to show the next
Theorem.{\parindent=0pt\v5

{\bf 4.1. Theorem} Suppose that the ideal $\I$ is generated by a
subset $\mathcal{G}\subset R_0\langle X\rangle$ which is a monic
Gr\"obner basis of $\I$ in $\RS$ with respect to a monomial ordering
$\prec$ on the standard $R$-basis $\B_R$ of $\RS$, where
$R_0\langle X\rangle =R_0\langle X_1,...,X_n\rangle$ is the free
$R_0$-algebra on $X_1,...,X_n$, and without loss of generality we
assume that $\LM (g)\ne 1$ for every $g\in\G$. Then the $R$-algebra
$A=\RS /\I =\RS /\langle\G\rangle$ can be endowed with the
$\Gamma$-filtration $FA$ as constructed in Theorem 3.1,  and if the
$\Gamma$-filtration $FR$ of $R$ is separated, then $FA$ is
separated, i.e., if $0\ne a\in A$, then  there is a
$\gamma\in\Gamma$ such that $a\in
F_{\gamma}A-F_{<\gamma}A,~\hbox{where}~F_{<\gamma}A=\cup_{\gamma
'<\gamma}F_{\gamma'}A.$ In particular, if $1\in F_0R-F_{<0}R$ then
$1\in F_0A-F_{<0}A$.\vskip 6pt

{\bf Proof} By the assumption on $\I$ and $\G$, Theorem 3.1 assures
the existence of the $\Gamma$-filtration $FA$ of $A$. Let $N(\G )$
be the set of normal monomials in $\B_R$ (mod $\G$). Then, by
Theorem 2.4, the canonical image $\OV{N(\G )}$ of $N(\G )$ in $A=\RS
/\I =\RS /\langle\G\rangle$ forms a free $R$-basis for $A$, and
moreover, $1\in \OV{N(\G )}$. If $0\ne a\in A$, then $a$ can be
written uniquely as $a=\sum_{j=1}^m\mu_j\bar w_j$ with $\mu_j\in R$,
$\bar w_j\in\OV{N(\G )}$. Since $FR$ is separated by the assumption,
there are $\gamma_1,\gamma_2,...,\gamma_m\in\Gamma$ such that
$\mu_j\in F_{\gamma_j}R-F_{<\gamma_j}R$, $1\le j\le m$. Assuming
$\gamma_1\ge\gamma_2\ge\cdots\ge\gamma_m$ in the totally ordered
monoid $\Gamma$, we have $\mu_j\in F_{\gamma_1}R$, $1\le j\le m$,
and hence $a\in F_{\gamma_1}A$. If there were some $\tau\in\Gamma$
with $\tau <\gamma_1$, such that $a\in F_{\tau}A$, then
$a=\sum^n_{i=1}\lambda_i\bar w_i$ with $\lambda_i\in F_{\tau}R$ and
$\bar w_i\in \OV{N(\G )}$. Comparing the coefficients on both sides
of the equality
$$\sum^m_{j=1}\mu_j\bar w_j=a=\sum^n_{i=1}\lambda_i\bar w_i,$$ we
would get $\mu_j\in F_{\tau}R$, $1\le j\le m$, in particular
$\mu_1\in F_{\tau}R$ with $\tau <\gamma_1$,  which is a
contradiction. Therefore $a\in F_{\gamma_1}A-F_{<\gamma_1}A$. This
shows that $FA$ is separated.}\par Finally, noticing $1\in\OV{N(\G
)}$, if $1\in F_0R-F_{<0}R$, then by the construction of $F_0A$ and
a similar argument as above we get $1\in F_0A-F_{<0}A$. \QED\v5

\section*{5. Realizing Good Reductions for $A$ by Gr\"obner Bases over $D\subset R$}
Let $R$ be an arbitrary commutative ring, $\RS =R\langle
X_1,...,X_n\rangle$ the free $R$-algebra on $X_1,...,X_n$, $\I$ an
ideal of $\RS$, and $A=\RS /\I$. In this section we generalize the
notion of a good reduction (in the sense of [MVO]) to the
$R$-algebra $A$, and we realize this property by using monic
Gr\"obner bases.\v5
\def\DX{D\langle X\rangle}\v5

Let $D$ be {\it any subring} of $R$ which has the same identity
element 1 as that of $R$, and let  $D\langle X\rangle =D\langle
X_1,...,X_n\rangle$ be the free $D$-algebra on $X_1,...,X_n$. In
what follows we use $B_R$, respectively $B_D$, to denote the
standard $R$-basis of $\RS$, respectively the standard $D$-basis of
$D\langle X\rangle$.  Considering the $D$-subalgebra
$$\Lambda = D\langle X\rangle +\I /\I$$ of $A$, we have $R\Lambda
=A$ and $\Lambda\cong D\langle X\rangle /\I\cap D\langle X\rangle$.
Observe that the exact sequence
$$0~\mapright{}{}~\I~\mapright{}{}~\RS~\mapright{\pi}{}~A~\mapright{}{}~0$$
and the canonical $D$-algebra epimorphism $\pi_D$: $D\langle
X\rangle\r \Lambda$ give rise to the exact sequence
$$0~\mapright{}{}~\I\cap D\langle X\rangle~\mapright{}{}~D\langle X\rangle~\mapright{\pi_D}{}~
\Lambda~\mapright{}{}~0$$ Let $\omega$ be {\it any  proper ideal} of
$D$ and $k=D/\omega$. Then the $k$-algebra epimorphism
$\pi_{\omega}$: $D\langle X\rangle /\omega D\langle
X\rangle\r\Lambda /\omega\Lambda$ induced by $\pi_D$ yields the
exact sequence
$$0~\mapright{}{}~\frac{\I\cap D\langle X\rangle+\omega D\langle X\rangle}{\omega D\langle X\rangle}
~\mapright{}{}~D\langle X\rangle /\omega D\langle
X\rangle~\mapright{\pi_{\omega}}{}~ \Lambda
/\omega\Lambda~\mapright{}{}~0$$ It is clear that if the ideal $\I$
of $\RS$ is generated by a subset $\G\subset\DX$ but $\G\not\subset
\omega\DX$, such that $\I\cap\DX=\langle\G\rangle$ as an ideal of
$\DX$, then $\Lambda /\omega\Lambda\cong k\langle X\rangle
/\langle\OV{\G}\rangle$ as $k$-algebras, where $\OV{\G}$ is the
canonical image of $\G$ in $\DX /\omega\DX$. {\parindent=0pt\v5

{\bf 5.1. Definition} (Compare with the definition given in [MVO]
Section 2) Let $\I =\langle \G\rangle$ be the ideal of  $\RS$
generated by a subset $\G\subset\DX$. If, as an ideal of $\DX$,
$\I\cap\DX=\langle\G\rangle$, then we say that the $D$-algebra
$\Lambda =\DX +\I /\I$ defines a {\it good reduction} for the
$R$-algebra $A=\RS /\I$.\v5

{\bf 5.2. Theorem} With the notation as before, if the ideal $\I$ is
generated by a subset $\mathcal{G}\subset D\langle X\rangle$ which
is a monic Gr\"obner basis of $\I$ in $\RS$ with respect to a
monomial ordering $\prec$ on the standard $R$-basis $\B_R$ of $\RS$,
then the following statements hold.\par

(i) $\G$ is a Gr\"obner basis for the ideal $\I \cap\DX$ in $\DX$
with respect to the same monomial ordering $\prec$ on the standard
$D$-basis $\B_D$ of $\DX$. Hence the ideal $\I \cap\DX$ of $\DX$ is
generated by $\G$, i.e.,  $\DX\cap\I =\langle\G\rangle$ holds in
$\DX$. Moreover, the set of normal monomials in $\B_D$ (mod $\G$) is
the same as the set of normal monomials in $\B_R$ (mod $\G$),
denoted $N(\G)$.

(ii) The $D$-algebra $\Lambda =\DX +\I /\I$ defines a good reduction
for  the $R$-algebra $A=\RS /\I $.\par

(iii) For any ideal $\omega$ of $D$ such that
$\G\not\subset\omega\DX$, we have $\Lambda /\omega\Lambda\cong
k\langle X\rangle /\langle\OV{\G}\rangle$ as $k$-algebras.\vskip 6pt

{\bf Proof} Although (i) is a consequence of Proposition 2.7, we
prefer giving a direct proof here. First note that $\B_R=\B_D$. If
$f\in\I \cap\DX$, say $f=\sum_id_iu_i$ with $d_i\in D$ and
$u_i\in\B_D$, then since $\G\subset\DX$ is a monic Gr\"obner basis
for the ideal $\I$ in $\KX$ with respect to the monomial ordering
$\prec$ on $\B_R$, we have $\OV{f}^{\G}=0$ in $\DX$ with respect to
the same monomial ordering $\prec$ on $\B_D$, that is, the division
of $f$ by $\G$ produces a Gr\"obner representation
$f=\sum_{i,j}d_{ij}w_{ij}g_jv_{ij}$ with $d_{ij}\in D$,
$w_{ij},v_{ij}\in\B_D =\B_R$, and $g_j\in\G$. Hence $\G$ is a
Gr\"obner basis for the ideal $\I \cap\DX$ in $\DX$ with respect to
the same monomial ordering $\prec$ on $\B_D$. It turns out that  $\I
\cap\DX =\langle\G\rangle$ in $\DX$, and that the set of normal
monomials in $\B_D$ (mod $\G$) is the same as the set of normal
monomials in $\B_R$ (mod $\G$) .}\par

(ii) and (iii) are clear enough by (i) and the discussion made
before Definition 5.1.\QED}\v5

\section*{6. Realizing Valuation Extensions of $A$ by Gr\"obner Bases over $\O_v$}
In this section, we apply the results of previous sections to
proving the main result of this paper.\v5

We first recall some basics on valuations, especially some
fundamental results concerning valuation extensions via
filtered-graded structures (cf. [Sc], [Coh], [LVO2], [MVO], [VO]).
Let $(\Gamma ,+;~<)$ be a totally ordered {\it abelian additive
group}  with the neutral element 0 and the total ordering $<$, and
let $A$ be an arbitrary (commutative or noncommutative) ring with 1.
A valuation $v$ of $A$ is a surjective function $A\r
\Gamma\cup\{\infty\}$, where the symbol $\infty$ plays
conventionally the role such that $\gamma<\infty$, $\infty +\gamma
=\gamma +\infty =\infty$ for every $\gamma\in\Gamma$, $\infty+\infty
=\infty$, and $\infty -\infty =0$, such that for $a,b\in A$, the
following conditions are satisfied:\par

(V1) $v(a)=\infty$ if and only if $a =0$;\par

(V2) $v(ab)=v(a)+v(b)$;\par

(V3) $v(a+b)\ge\min(v(a),v(b))$.{\parindent=0pt\par

If $A$ is an Ore domain with a valuation function $v$ as above, then
$v$ can be extended in a unique way to the (skew-)field $\Delta$ of
fractions of $A$ subject to the rule:
$$v(ab^{-1})=v(a)-v(b)~\hbox{for}~a,b\in A~\hbox{with}~b\ne 0.$$}\par

Valuation theory is closely related to filtered-graded structures.
If a ring $A$ has a valuation $v$: $A\r\Gamma\cup\{\infty\}$, then
$v$ determines an exhaustive $\Gamma$-filtration $F^vA=\{
F_{\gamma}^vA\}_{\gamma\in\Gamma}$ for $A$ by putting
$$F_{\gamma}^vA=\{ a\in A~|~v(a)\ge -\gamma\},~\gamma\in\Gamma ,$$ i.e.,
$F^vA$ satisfies (F1) -- (F4) as described in the beginning of
Section 3. For the convenience of later use we also note three more
properties of $F^vA$ as follows.{\parindent=0pt \v5

{\bf 6.1. Proposition} The exhaustive $\Gamma$-filtration $F^vA$ of
$A$ defined above has the following properties: \par

(i)  For  $0\ne a\in A$, $v(a)=-\gamma$ if and only if $a\in
F^v_{\gamma}A-F^v_{<\gamma}A$, where $F^v_{<\gamma}A=\cup_{\gamma
'<\gamma}F^v_{\gamma '}A$. Hence $F^vA$ is separated. In particular,
$1\in F_0A-F_{<0}A$.\par

(ii) If $a\in F^v_{\gamma}A-F^v_{<\gamma}A$ and $a$ is invertible,
then $a^{-1}\in F_{-\gamma}^vA-F_{<-\gamma}A$.\par

(iii) Let $G(A)=\oplus_{\gamma\in\Gamma}G(A)_{\gamma}$ be the
associated $\Gamma$-graded ring determined by $FA$, where
$G(A)_{\gamma}=F_{\gamma}A/F_{<\gamma}A$ for every
$\gamma\in\Gamma$. Then $G(A)$ is a domain and thereby $A$ is a
domain. \par\QED\v5

Conversely, let $A$ be a $\Gamma$-filtered ring with an exhaustive
$\Gamma$-filtration $FA=\{ F_{\gamma}A\}_{\gamma\in\Gamma}$ in the
sense of Section 3. If $FA$ is separated, i.e., $0\ne a\in A$
implies that there is a $\gamma\in\Gamma$ such that $a\in
F_{\gamma}A-F_{<\gamma}A$, where $F_{<\gamma}A=\cup_{\gamma
'<\gamma}F_{\gamma '}A$, in particular, we insist  that $1\in
F_0A-F_{<0}A$. then the {\it degree function} on $A$ can be defined
by setting
$$\begin{array}{cccl} d:&A&\mapright{}{}&\Gamma\cup\{\infty\}\\
&a&\mapsto&\left\{\begin{array}{ll} \gamma ,&\hbox{if}~a\in
F_{\gamma}A-F_{<\gamma}A,\\
-\infty ,&\hbox{if}~a=0\end{array}\right.\end{array}$$ Furthermore,
consider the associated $\Gamma$-graded ring
$G(A)=\oplus_{\gamma\in\Gamma}G(A)_{\gamma}$ of $A$ determined by
$FA$, where $G(A)_{\gamma}=F_{\gamma}A/F_{<\gamma}A$ for every
$\gamma\in\Gamma$. If $G(A)$ is a domain (hence $A$ is a domain),
then the function defined by setting
$$\begin{array}{cccl} v:&A&\mapright{}{}&\Gamma\cup\{\infty\}\\
&a&\mapsto&-d(a)\end{array}$$ is a valuation function on $A$.}\par

In conclusion, the next theorem summarizes the principle of
valuation extensions via filtered-graded structures.
{\parindent=0pt\v5

{\bf 6.2. Theorem} Let $A$ be a $\Gamma$-filtered ring with an
exhaustive and separated filtration $FA=\{
F_{\gamma}A\}_{\gamma\in\Gamma}$ such that $1\in F_0A-F_{<0}A$. Then
the following statements hold. \par

(i) The degree function $d(x)$ on $A$ defines a valuation function
$v(x)=-d(x)$ on $A$ if and only if $G(A)$ is a domain.\par

(ii) Suppose that $G(A)$ is a domain (hence $A$ is a domain) and
the (skew-)field $\Delta$ of fractions of $A$ exists, then $v$ can
be uniquely extended to a valuation function on $\Delta$, or
equivalently, the $\Gamma$-filtration $FA$ can be extended to an
exhaustive and separated $\Gamma$-filtration $F\Delta$ such that
$F_{\gamma}A=A\cap F_{\gamma}\Delta$ for every $\gamma\in\Gamma$,
and moreover, $G(\Delta )$ is a $\Gamma$-graded (skew-)field in the
sense that every nonzero homogenous element of $G(\Delta )$ is
invertible. \par\QED}\v5

Now, let $K$ be a field  $v$: $K\r \Gamma\cup\{\infty\}$ be a
valuation of $K$. Then, by the foregoing discussion, $v$ determines
an exhaustive and separated $\Gamma$-filtration $F^vK=\{
F^v_{\gamma}K\}_{\gamma\in\Gamma}$ with $F_{\gamma}^vK=\{ \lambda\in
K~|~v(\lambda )\ge -\gamma\}$ for $K$, which has all properties as
described in Proposition 6.1. Moreover, since $K$ is a field, it
follows from Proposition 6.1 that $F^vK$ is a strong
$\Gamma$-filtration in the sense of [LVO1], i.e.,
$F_{\gamma_1}^vK\cdot F_{\gamma_2}^vK=$ $F^v_{\gamma_1+\gamma_2}K$,
in particular, $F^v_{\gamma}K\cdot F^v_{-\gamma}K=F^v_0K$,  for all
$\gamma_1,\gamma_2,\gamma\in\Gamma$, and this leads to the fact that
the associated graded ring  $G(K)$ of $K$ is a strongly
$\Gamma$-graded ring in the sense of [NVO], i.e.,
$G(K)_{\gamma_1}\cdot G(K)_{\gamma_2}=G(K)_{\gamma_1+\gamma_2}$, in
particular, $G(K)_{\gamma}\cdot G(K)_{-\gamma}=G(K)_0$, for all
$\gamma_1,\gamma_2,\gamma\in\Gamma$. Noting that the valuation ring
of $K$ associated to $v$ is local ring $\O_v=\{ \lambda\in
K~|~v(\lambda )\ge 0\}$ with the unique maximal ideal $\M_v=\{
\lambda\in K~|~v(\lambda )> 0\}$, by the definition of $F^vK$ we
have $F^v_0K=\O_v$ and $F^v_{<0}K=\M_v$. Thus, $G(K)_0=k=\O_v /\M_v$
is the residue field of $\O_v$ and $G(K)$ is indeed a commutative
$\Gamma$-graded field in the sense that every nonzero homogeneous
element of $G(K)$ is invertible.  \par

Next, consider the free $K$-algebra $\KS =K\langle
X_1,...,X_n\rangle$ on $X_1,...,X_n$. Let $\I$ be an ideal of $\KS$,
$A=\KS /\I$, and $\Lambda =\O_v\langle X\rangle +\I /\I$, where
$\O_v\langle X\rangle =\O_v\langle X_1,...,X_n\rangle$ is the free
$\O_v$-algebra on $X_1,...,X_n$. If $\G$ is a Gr\"obner basis for
$\I$ in $\KS$ with respect to a monomial ordering $\prec$ on the
standard $K$-basis $\B_K$ of $\KS$, as before we write $N(\G )$ for
the set of normal monomials in $\B_K$ (mod $\G$), and we write
$\OV{N(\G )}$ for the canonical image of $N(\G )$ in $A$, which is
known a $K$-basis for $A$ (Theorem 2.4). {\parindent=0pt\v5

{\bf 6.3. Theorem} With all notations as we fixed so far, suppose
that the ideal $\I$ is generated in $\KS$ by a monic Gr\"obner basis
$\G\subset \O_v\langle X\rangle$ with respect to a monomial ordering
$\prec$ on $\B_K$, and without loss of generality we assume that
$\LM (g)\ne 1$ for every $g\in\G$. Then the following statements
hold.\par

(i) The $K$-algebra $A=\KS /\I=\KS /\langle\G\rangle$ can be endowed
with the exhaustive $\Gamma$-filtration  $F^vA=\{
F^v_{\gamma}A\}_{\gamma\in\Gamma}$ by putting
$$F^v_{\gamma}A=\left\{\left. a=\sum_i\lambda_i\bar w_i~\right |~\lambda_i\in F^v_{\gamma}K,~
\bar w_i\in \OV{N(\G)}\right\} , ~\gamma\in\Gamma ,$$ such that
$F_{\gamma}K=K\cap F_{\gamma}A,~\gamma\in\Gamma,$ that is, $F^vK$
extends naturally to $F^vA$.\par

(ii) The $\Gamma$-filtration $F^vA$ obtained in (i) is separated,
i.e., if $0\ne a\in A$, then  there is a $\gamma\in\Gamma$ such that
$a\in F^v_{\gamma}A-F^v_{<\gamma}A,$ in particular $1\in
F^v_0A-F^v_{<0}A$, where $F^v_{<\gamma}A=\cup_{\gamma
'<\gamma}F^v_{\gamma'}A.$\par

(iii) The $\Gamma$-filtration $F^vA$ obtained in (i) has
$$\begin{array}{l} F_0^vA=\O_v\langle X\rangle /\langle\G\rangle\cong \O_v\langle X\rangle
+\I /\I =\Lambda,\\ F^v_{<0}A=\M_vF^v_0A=\M_v\Lambda.\end{array}$$
Hence the associated $\Gamma$-graded $K$-algebra
$G(A)=\oplus_{\gamma\in\Gamma}G(A)_{\gamma}$ of $A$ determined by
$F^vA$ has $G(A)_0=F^v_0A/F^v_{<0}A=\Lambda /\M_v\Lambda$\par

(iv) The $\Gamma$-filtration $F^vA$ obtained in (i) is a strong
filtration, i.e., $F_{\gamma_1}^vA\cdot
F^v_{\gamma_2}A=F^v_{\gamma_1+\gamma_2}A$ for all
$\gamma_1,\gamma_2\in \Gamma$, and the associated $\Gamma$-graded
$K$-algebra $G(A)=\oplus_{\gamma\in\Gamma}G(A)_{\gamma}$ with
$G(A)_{\gamma}=F^v_{\gamma}A/F^v_{<\gamma}A$ is strongly
$\Gamma$-graded, i.e., $G(A)_{\gamma_1}\cdot
G(A)_{\gamma_2}=G(A)_{\gamma_1+\gamma_2}$ for all
$\gamma_1,\gamma_2\in\Gamma$.\par

(v) $\G$ is a Gr\"obner basis for the ideal $\I \cap\DX$ of $\DX$
with respect to the same monomial ordering $\prec$ on the standard
$\O_v$-basis $\B_{\O_v}$ of $\O_v\langle X\rangle$. It follows that
$\DX\cap\I =\langle\G\rangle$ holds in $\DX$, and thereby the
$\O_v$-algebra $\Lambda =\O_v\langle X\rangle +\I /\I$ defines a
good reduction for  the $K$-algebra $A=\KS /\I $ in the sense of
Definition 5.1.\par

(vi) If $\G\not\subset\M_v\O_v\langle X\rangle$, then $\Lambda
/\M_v\Lambda\cong k\langle X\rangle /\langle\OV{\G}\rangle$ as
$k$-algebras, where $k=\O_v/\M_v$ is the residue field of $\O_v$,
$k\langle X\rangle =k\langle X_1,...,X_n\rangle$ is the free
$k$-algebra on $X_1,...,X_n$,  and $\OV{\G}$ is the canonical image
of $\G$ in $\O_v\langle X\rangle /\M_v\O_v\langle X\rangle$\par

(vii) If $\G\not\subset\M_v\O_v\langle X\rangle$ and $k\langle
X\rangle /\langle\OV{\G}\rangle$ is a domain, then $G(A)$ is a
domain and thereby $A$ is a domain. It follows that $F^vA$
determines a valuation function $A\rightarrow \Gamma\cup\{\infty\}$,
and thereby  $v$ extends naturally to a valuation function on the
(skew-)field $\Delta$ of fractions of $A$ provided $\Delta$ exists.
\vskip 6pt

{\bf Proof} Note that the $\Gamma$-filtration $F^vK$ of $K$
determined by the valuation $v$: $K\r \Gamma\cup\{\infty\}$ is
exhaustive and separated. Moreover, $F^v_0K=\O_v$ and
$F^v_{<0}K=\M_v$.} \par

(i) and (ii) follow from Theorem 3.1, Theorem 4.1, and Proposition
6.1.\par

(iv) follows from (ii), (iii), and Proposition 6.1.\par

(iii), (v), and (vi) follow from Theorem 5.2.\par

(vii) By the foregoing (iii), (vi) and (iv), $G(A)$ is now a
strongly $\Gamma$-graded algebra with $G(A)_0=\Lambda
/\M_v\Lambda\cong k\langle X\rangle /\langle\OV{\G}\rangle$.  If
$k\langle X\rangle /\langle\OV{\G}\rangle$ is a domain, then $G(A)$
is a domain and thereby $A$ is a domain. It follows from Theorem
6.2. that the last assertion holds.\QED}\v5

Let $K[x_1,...,x_n]$ be the commutative polynomial $K$-algebra in
$n$ variables over a field $K$. Noticing $K[x_1,...,x_n]\cong
K\langle X_1,...,X_n\rangle /\langle\G\rangle$ with $\G=\{
X_jX_i-X_iX_j~|~1\le i<j\le n\}$ a Gr\"obner basis for the ideal
$\langle \G\rangle$, Theorem 6.3 has an immediate application to
$K[x_1,...,x_n]$. {\parindent=0pt\v5

{\bf 6.4. Corollary} Let $K$ be a field. Then every valuation $v$ on
$K$ extends naturally to a valuation function on $K[x_1,...,x_n]$
and further to a valuation function on the field of rational
functions $K(x_1,...,x_n)$.
\par\QED}\v5

More generally, as it was pointed out in ([Li3], Section 1,
Remark(iv)), Proposition 2.7 of previous Section 2 is valid for
getting monic Gr¡§obner bases in a commutative polynomial ring
$R[x_1,...,x_n]$ over an arbitrary commutative ring $R$ when overlap
elements are replaced by S-polynomials. It follows that the results
of Sections 3 -- 5 and Theorem 6.3 are also valid for commutative
algebras over a field $K$ after replacing $\KS$ by $K[x_1,...,x_n]$.
For instance, let $\O_v$ be a valuation ring of $K$ associated to a
valuation $v$ of $K$, and let $A=K[x_1,...,x_n]/I$ be the coordinate
ring of an affine variety $V(I)\subset K^n$. If the ideal $I$ is
generated by a subset $\G\subset\O_v[x_1,...,x_n]$ which is a
Gr\"obner basis of $I$ in $K[x_1,...,x_n]$ with respect to a
monomial ordering on $K[x_1,...,x_n]$, then Theorem 6.3 holds for
$A$ after replacing $\KS$ by $K[x_1,...,x_n]$.\par

\parindent=0pt\vskip 1truecm

\centerline{References}\parindent=1truecm\par

\item{[BVO]} C. Baetica and F. Van Oystaeyen, Valuation extensions of filtered and graded algebras,
{\it Comm. Alg.}, 3(34)(2006), 829--840.

\item{[Coh]} P.M. Cohn, {\it Algebra} I \& II, Hohn Wiley and Sons Ltd., 1982.

\item{[G-I]} T. Gateva-Ivanova, Binomial skew polynomial rings, Artin-Schelter
regularity, and binomial solutions of the Yang-Baxter equation, {\it
Serdica Math. J.}, 30(2004), 431--470. arXiv:0909.4707.

\item{[Gr]} E. L. Green, Noncommutative Gr¡§obner bases and projective
resolutions, in: {\it Proceedings of the Euroconference
Computational Methods for Representations of Groups and Algebras},
Essen, 1997, (Michler, Schneider, eds), Progress in Mathematics,
Vol. 173, Basel, Birkha¡§user Verlag, 1999, 29--60.

\item{[Laf]} G. Laffaille, Quantum binomial algebras, Colloquium on Homology and
Representation Theory (Spanish) (Vaquerfas, 1998). Bol. Acad. Nac.
Cienc. (C\'ordoba) 65 (2000), 177--182.

\item{[Li1]} H. Li, A note on the extension of discrete valuations to affine domains, {\it Comm. Alg.},
25(1997), 1805--1816.

\item{[Li2]} H. Li, {\it Noncommutative Gr\"obner Bases and
Filtered-Graded Transfer}, LNM, 1795, Springer-Verlag, 2002.

\item{[Li3]} H. Li,  Algebras defined by monic Gr\"obner bases over rings,
arXiv:0906.4396.

\item{[LVO1]} H. Li and F. Van Oystaeyen, Strongly filtered rings applied to Gabber's integrability
theorem and modules with regular singularities, in: {\it Proc. Sem.
Malliavin}, LNM., 1404, Springer-Verlag, 1988, 296--322

\item{[LVO2]} H. Li and F. Van Oystaeyen, Filtration on simple artinian rings, {\it J. Alg.}, 132(1990), 361--376.

\item{[LVO3]} H. Li and F. Van Oystaeyen, Reductions and global dimension of quantized algebras over a
regular commutative domain, {\it Comm. Alg.}, 26(4)(1998),
1117--1124.

\item{[Mor]} T. Mora, An introduction to commutative and
noncommutative Gr\"obner bases, {\it Theoretic Computer Science},
134(1994), 131--173.

\item{[MVO]} M. Hussein and F. Van Oystaeyen, Discrete valuations extend to certain algebras of quantum type,
{\it Comm. Alg.}, 24(8)(1996), 2551--2566.

\item{[NVO]} C. N$\check{\rm a}$st$\check{\rm a}$sescu, F. Van Oystaeyen, {\it Graded ring theory}, Math. Library 108, North Holland,
Amsterdam,1982.

\item{[Sc]} O.F.G. Schilling, {\it The theory of valuations}, Mathematical Surveys, AMS, 1950.

\item{[VO]} F. Van Oystaeyen, {\it Algebraic geometry for associative algebras}, Monographs and Textbooks in Pure and Applied
Mathematics, Vol. 232, Marcel Dekker, 2000.

\item{[VOW1]} F. Van Oystaeyen and L. Willaert, Grothendieck topology,
coherent sheaves and Serre's theorem for schematic algebras,    {\it
J. Pure App. Algebra}, 1(104)(1995), 109-122.

\item{[VOW2]} F. Van Oystaeyen and L. Willaert, Valuations on extensions of Weyl skewfields,
{\it J. Alg.}, 183(1996), 359--364.

\item{[VOW3]} F. Van Oystaeyen and L. Willaert, Examples and quantum sections of schematic algebras,
{\it J. Pure App. Algebra},  2(120)1997, 195-211.

\end{document}